\documentclass[10pt]{article}
\usepackage{amsfonts}
\usepackage{amssymb}
\usepackage{amsthm}
\usepackage{amsmath}
\usepackage{latexsym}

\newtheorem{theorem}{Theorem}[section]
\newtheorem{lemma}{Lemma}[section]
\newtheorem{corollary}{Corollary}[section]

\numberwithin{equation}{section}
\begin{document}

\title{Asymptotics for a Discrete-time Risk Model with the Emphasis on Financial Risk}
\author{Enkelejd Hashorva\thanks{%
Department of Actuarial Science, Faculty of Business and Economics,
University of Lausanne, UNIL-Dorigny 1015 Lausanne, Switzerland. Email:
Enkelejd.Hashorva@unil.ch} \ and Jinzhu Li\thanks{%
School of Mathematical Science and LPMC, Nankai University, Tianjin 300071,
P.R. China}}
\date{{\small \today}}
\maketitle

\begin{abstract}
This paper focuses on a discrete-time risk model in which both insurance
risk and financial risk are taken into account. We study the asymptotic
behaviour of the ruin probability and the tail probability of the aggregate
risk amount. Precise asymptotic formulas are derived under weak moment
conditions on involved risks. The main novelty of our results lies in the
quantification of the impact of the financial risk.


\end{abstract}

\baselineskip15pt

\section{Introduction and Preliminaries}

In this paper, for every $i\geq 1$, let $X_{i}$ be an insurer's net loss
(the total amount of claims less premiums) within period $i$ and let $Y_{i}$
be the stochastic discount factor (the reciprocal of the stochastic return
rate) over the same time period. Then the stochastic present values of
aggregate net losses of the insurer can be specified as 
\begin{equation}
S_{0}=0,\qquad S_{n}=\sum_{i=1}^{n}X_{i}\prod_{j=1}^{i}Y_{j},\qquad n\geq 1,
\label{S}
\end{equation}%
with their maxima 
\begin{equation}
M_{n}=\max_{0\leq k\leq n}S_{k},\qquad n\geq 1.  \label{M}
\end{equation}%
We are concerned with the asymptotic behaviour of the tail probabilities $%
\mathbb{P}\left( S_{n}>x\right) $ and $\mathbb{P}\left( M_{n}>x\right) $ as $%
x\rightarrow \infty $, in which $\mathbb{P}\left( M_{n}>x\right) $ coincides
with the insurer's finite-time ruin probability within period $n$ given that
the initial wealth is $x$.

In the literature $\{X_{i};i\geq 1\}$ and $\{Y_{i};i\geq 1\}$\ are usually
called the insurance risk and the financial risk, respectively. Under
certain independence or identical distribution assumptions imposed on $X_{i}$%
's and $Y_{i}$'s, the asymptotic tail behaviour of $S_{n}$ and $M_{n}$ has
been extensively studied by many researchers. See, e.g., Tang and
Tsitsiashvili (2003, 2004), Konstantinides and Mikosch (2005), Tang (2006),
Zhang et al. (2009), Chen (2011), and Yang and Wang (2013) for some recent
findings. Since the products of $Y_{i}$'s appearing in (\ref{S}) essentially
cause technical problems in the derivation of explicit asymptotic formulas,
most of existing works assumed that the financial risk is dominated by the
insurance risk, i.e., the tails of $Y_{i}$'s are lighter than the tails of $%
X_{i}$'s, usually through imposing strong moment conditions on $Y_{i}$'s.
Then the problem becomes relatively tractable and the final results are
mainly determined by the tails of $X_{i}$'s.

However, as shown by empirical data and the most recent financial crisis,
the financial risk may impair the insurer's solvency as seriously as does
the insurance risk and, hence, it should not be underestimated as before;
see Norberg (1999), Frolova et al. (2002), Kalashnikov and Norberg (2002),
and Pergamenshchikov and Zeitouny (2006). Therefore, in the current
contribution, we focus on the other directions where the financial risk
dominates the insurance risk or no dominating relationship exists between
the two kinds of risk. We aim at capturing the impact of the financial risk
(the products of $Y_{i}$'s) on the tail behaviour of $S_{n}$ and $M_{n}$.
Loosening some independence and identical distribution constraints, we
derive precise asymptotic formulas under weak moment conditions on $Y_{i}$'s
and $X_{i}$'s.

Throughout this paper, an underlying assumption is the following:

\vskip0.2cm

\noindent \textbf{Assumption A.} $\{X_{i};i\geq 1\}$ is a sequence of
real-valued rv's (random variables) with distribution functions $F_{i}$'s, $%
\{Y_{i};i\geq 1\}$ is a sequence of positive and independent rv's with
distribution functions $G_{i}$'s, and $\{X_{i};i\geq 1\}$ and $\{Y_{i};i\geq
1\}$ are mutually independent.

\vskip0.2cm

It is worth mentioning that, if we further assume that both $\{X_{i};i\geq
1\}$ and $\{Y_{i};i\geq 1\}$ are sequences of iid (independent and
identically distributed) rv's in (\ref{S}), then there is a natural
connection between this discrete-time risk model and the general bivariate L%
\'{e}vy-driven risk model with the form%
\begin{equation*}
U_{t}=\int_{0}^{t}\mathrm{e}^{Q_{s}}\mathrm{d}P_{s},\qquad t\geq 0,
\end{equation*}%
where $\{Q_{s};s\geq 0\}$ and $\{P_{s};s\geq 0\}$ are two independent L\'{e}%
vy processes; see Paulsen (1993, 2008), Hao and Tang (2012), and the
references therein. To see this, arbitrarily embed an increasing sequence of
stopping times, say $\{\tau _{i};i\geq 1\}$, to the continuous-time model.
Then, after such a discretization procedure, $U_{\tau _{n}}$ takes the form
as $S_{n}$ in (\ref{S}). Due to this reason, the results obtained in this
paper can provide us with some valuable insights to the general bivariate L%
\'{e}vy-driven case.

We restrict our discussions within the scope that $Y_{i}$'s are regularly
varying. A real-valued rv $Z$ with distribution function $H$ is said to be
regularly varying if its survival function $\overline{H}=1-H$ is regularly
varying at infinity, i.e., $\lim_{x\rightarrow \infty }\overline{H}(xy)/%
\overline{H}(x)=y^{-\alpha }$ for every $y>0$ and some $\alpha \geq 0$. In
this case, we write $Z\in \mathcal{R}_{-\alpha }$ or $\overline{H}\in 
\mathcal{R}_{-\alpha }$. A positive function regularly varying with $\alpha
=0$ is also called slowly varying function.\ See Bingham et al.\ (1987),
Resnick (1987), or Embrechts et al.\ (1997) for more details on regularly
varying functions.

Hereafter, all limit relations hold as $x\rightarrow \infty $ unless
otherwise specified. For two positive functions $a(\cdot )$ and $b(\cdot )$,
we write $a(x)\gtrsim b(x)$ or $b(x)\lesssim a(x)$ if $\liminf_{x\rightarrow
\infty }a(x)/b(x)\geq 1$ and write $a(x)\sim b(x)$\ if both $a(x)\lesssim
b(x)$ and $a(x)\gtrsim b(x)$.

Our first result below shows that, in a special case of regular variation,
the moment conditions of involved rv's can be dropped thanks to a Rootz\'{e}%
n-type lemma stated in Section 3 (Lemma \ref{product}).

\begin{theorem}
\label{main3}Under Assumption $\mathbf{A}$, let $X_{i}$'s be independent.
If, for every $i\geq 1$, $\overline{F}_{i}(x)\sim \ell _{i}^{\ast }(\ln
x)\cdot \left( \ln x\right) ^{\gamma ^{\ast }-1}x^{-\alpha }$ and $\overline{%
G}_{i}(x)\sim \ell _{i}(\ln x)\left( \ln x\right) ^{\gamma _{i}-1}x^{-\alpha
}$ for some positive constants $\alpha ,\gamma ^{\ast },\gamma _{i}$ and
some slowly varying functions $\ell _{i}^{\ast }(\cdot ),\ell _{i}(\cdot )$
then, for every $n\geq 1$, letting $\bar{\gamma}_{n}=\gamma ^{\ast
}+\sum_{i=1}^{n}\gamma _{i}$, we have 
\begin{equation}
\mathbb{P\!}\left( S_{n}\mathbb{\!}>\mathbb{\!}x\right) \mathbb{\!\!}\sim 
\mathbb{\!P\!}\left( M_{n}\mathbb{\!}>\mathbb{\!}x\right) \mathbb{\!\!}\sim 
\mathbb{\!P\!\!}\left( \mathbb{\!}X_{n}\mathbb{\!}\prod_{j=1}^{n}\mathbb{\!}%
Y_{j}\mathbb{\!}>\mathbb{\!}x\mathbb{\!}\right) \mathbb{\!\!}\sim \mathbb{\!}%
\frac{\alpha ^{n}\Gamma (\gamma ^{\ast })\prod_{i=1}^{n}\Gamma (\gamma _{i})%
}{\Gamma \left( \bar{\gamma}_{n}\right) }\ell _{n}^{\ast }(\ln x)\mathbb{\!}%
\left( \prod\limits_{i=1}^{n}\ell _{i}(\ln x)\right) \mathbb{\!}(\ln x)^{%
\bar{\gamma}_{n}-1}x^{-\alpha }.  \label{Th5}
\end{equation}
\end{theorem}

\noindent \textbf{Remark 1.1.} A well-known folklore in risk theory is that
the ruin of an insurer, i.e., the tail of $M_{n}$, will be determined by one
of the insurance risk and the financial risk which has a heavier tail.
Nevertheless, Theorem \ref{main3} provides a counterexample violating the
folklore. To see this more clearly, let both $\{X_{i};i\geq 1\}$ and $%
\{Y_{i};i\geq 1\}$\ be sequences of iid rv's with common survival functions $%
\overline{F}(x)\sim \ell ^{\ast }(\ln x)\left( \ln x\right) ^{\gamma ^{\ast
}-1}x^{-\alpha }$ and $\overline{G}(x)\sim \ell (\ln x)\left( \ln x\right)
^{\gamma -1}x^{-\alpha }$, respectively. Then, according to the different
selections of $\gamma ^{\ast },\ell ^{\ast }(\cdot )$ and $\gamma ,\ell
(\cdot )$, Theorem \ref{main3} covers various asymptotic relationships
between $\overline{F}$ and $\overline{G}$. However, we have the unified
asymptotic expansion determined by both $\overline{F}$ and $\overline{G}$.

\vskip0.2cm

\noindent \textbf{Remark 1.2. }Tang and Tsitsiashvili (2003) gave a similar
result for $M_{n}$ in their Theorem 6.2. Their result does not cover, and is
not covered by, our Theorem \ref{main3}, since their conditions on $X_{i}$'s
and ours are mutually exclusive. However, their assumptions imply $\overline{%
F}(x)=o\left( \overline{G}(x)\right) $, whereas our Theorem \ref{main3}, as
stated in Remark 1.1, is valid for various relationships between $\overline{F%
}$ and $\overline{G}$.

\vskip0.2cm

Theorem \ref{main3} presents an elegant result which is due to the special
forms of $\overline{F}_{i}$'s and $\overline{G}_{i}$'s. In the subsequent
sections we focus on asymptotic analysis of $S_{n}$ and $M_{n}$ for general
regularly varying conditions, while the price to pay for it is the lack of
elegance and the high technicalities of the proofs. Our main results
presented in Theorem \ref{main2} below show that, as expected, similarly to
Theorem \ref{main3}, both $S_{n}$ and $M_{n}$ are regularly varying rv's
under some general conditions. Furthermore, we derive precise tail
asymptotics for both $S_{n}$ and $M_{n}$. One remarkable feature of our
Theorem \ref{main2} is the weakening of the moment assumptions commonly
imposed on $X_{i}$'s and $Y_{i}$'s in the literature.

The rest of the paper is organized as follows. Section 2 shows our main
theorem with several interesting remarks. Section 3 gives the lemmas and
proofs related to the results presented in Sections 1 and 2. As an appendix,
Section 4 discusses the constant weighted sums of the products of $Y_{i}$'s (%
$X_{i}\equiv c_{i}>0$ for every $i\geq 1$ in (\ref{S})), which model the
stochastic present values of some risk-free bond with fixed income $c_{i}$
in period $i$.\ We derive an asymptotic formula with the uniformity of the
constant weights in this case.

\section{Main Results and Remarks}

Hereafter, the summation and the product over an empty set of indices are
considered as $0$ and $1$, respectively. Moreover, to avoid triviality,
every individual real-valued rv is assumed to be not only concentrated on $%
(-\infty ,0]$. For a real number $a$, we write $a_{+}=a\vee 0$.

Under the framework specified in Assumption $\mathbf{A}$, we continue to
study the tail behaviour of $S_{n}$ and $M_{n}$ defined in (\ref{S}) and (%
\ref{M}). For the conciseness in writing and presentation, we further define%
\begin{equation*}
S_{0}^{(l)}=0,\qquad
S_{n}^{(l)}=\sum_{i=l}^{n+l-1}X_{i}\prod_{j=l}^{i}Y_{j},\qquad
l,n=1,2,\ldots ,
\end{equation*}%
and%
\begin{equation*}
M_{n}^{(l)}=\max_{0\leq k\leq n}S_{k}^{(l)},\qquad l,n=1,2,\ldots .
\end{equation*}%
Clearly, $S_{n}^{(l)}$ describes the stochastic present value at time $l-1$
of aggregate net losses occurring from time $l$ to time $n+l-1$. Note in
passing that $S_{n}^{(1)}=S_{n}$, $M_{n}^{(1)}=M_{n}$, and further 
\begin{equation}
S_{n}^{(l)}=Y_{l}\left( X_{l}+S_{n-1}^{(l+1)}\right) \text{\ and }%
M_{n}^{(l)}=Y_{l}\left( X_{l}+M_{n-1}^{(l+1)}\right) _{+},\qquad
l,n=1,2,\ldots .  \label{S-M}
\end{equation}

Our main results are given in the following Theorem \ref{main2}, in which
assertion (i) is valid for arbitrarily dependent $X_{i}$'s, assertion (ii)
drops the dominating relationship between $\overline{F}_{i}$'s and $%
\overline{G}_{i}$'s, and neither assertion (i) nor (ii) requires $\mathbb{E}%
\left( X_{i}\right) _{+}^{\beta }<\infty $ or $\mathbb{E}Y_{i}^{\beta
}<\infty $ for every $i\geq 1$ and some $\beta >\alpha $.

\begin{theorem}
\label{main2}Under Assumption $\mathbf{A}$, assume that $\overline{G}_{i}\in 
\mathcal{R}_{-\alpha }$ for every $i\geq 1$ and some $\alpha \geq 0$, and $%
\mathbb{E}Y_{i}^{\alpha }<\infty $ for every $i\geq 2$.

(i) If $X_{i}Y_{i}\in \mathcal{R}_{-\alpha }$ and 
\begin{equation}
\mathbb{P}\left( \left\vert X_{i}\right\vert >x\right) =o\left( \overline{G}%
_{i+1}(x)\right)  \label{assume}
\end{equation}%
for every $i\geq 1$ then, for every $n\geq 1$, $S_{n}\in \mathcal{R}%
_{-\alpha }$, $M_{n}\in \mathcal{R}_{-\alpha }$, and further 
\begin{equation}
\mathbb{P}\left( S_{n}>x\right) \sim \sum_{i=1}^{n-1}B_{n,i}\mathbb{P}\left(
\prod\limits_{j=1}^{i}Y_{j}>x\right) +\mathbb{P}\left(
X_{n}\prod\limits_{j=1}^{n}Y_{j}>x\right)  \label{Th1}
\end{equation}%
and 
\begin{equation}
\mathbb{P}\left( M_{n}>x\right) \sim \sum_{i=1}^{n-1}D_{n,i}\mathbb{P}\left(
\prod\limits_{j=1}^{i}Y_{j}>x\right) +\mathbb{P}\left(
X_{n}\prod\limits_{j=1}^{n}Y_{j}>x\right) ,  \label{Th2}
\end{equation}%
where%
\begin{equation*}
B_{n,i}=\mathbb{E}\left( X_{i}+S_{n-i}^{(i+1)}\right) _{+}^{\alpha }-\mathbb{%
E}\left( S_{n-i}^{(i+1)}\right) _{+}^{\alpha }\text{ and}\quad D_{n,i}=%
\mathbb{E}\left( X_{i}+M_{n-i}^{(i+1)}\right) _{+}^{\alpha }-\mathbb{E}%
\left( M_{n-i}^{(i+1)}\right) ^{\alpha }.
\end{equation*}

(ii) If $X_{i}$'s are independent and $\overline{F}_{i}\in \mathcal{R}%
_{-\alpha }$ with $\mathbb{E}\left( X_{i}\right) _{+}^{\alpha }<\infty $ for
every $i\geq 1$ then, for every $n\geq 1$, $S_{n}\in \mathcal{R}_{-\alpha }$%
, $M_{n}\in \mathcal{R}_{-\alpha }$, and further%
\begin{equation}
\mathbb{P}\left( S_{n}>x\right) \sim \sum_{i=1}^{n-1}\left( B_{n,i}-\mathbb{E%
}\left( X_{i}\right) _{+}^{\alpha }\right) \mathbb{P}\left(
\prod\limits_{j=1}^{i}Y_{j}>x\right) +\sum_{i=1}^{n}\mathbb{P}\left(
X_{i}\prod\limits_{j=1}^{i}Y_{j}>x\right)  \label{Th3}
\end{equation}%
and%
\begin{equation}
\mathbb{P}\left( M_{n}>x\right) \sim \sum_{i=1}^{n-1}\left( D_{n,i}-\mathbb{E%
}\left( X_{i}\right) _{+}^{\alpha }\right) \mathbb{P}\left(
\prod\limits_{j=1}^{i}Y_{j}>x\right) +\sum_{i=1}^{n}\mathbb{P}\left(
X_{i}\prod\limits_{j=1}^{i}Y_{j}>x\right) .  \label{Th4}
\end{equation}
\end{theorem}

One important theoretical merit of Theorem \ref{main2} lies in that, through
the transparent expansions (\ref{Th1})--(\ref{Th4}), it gives new criteria
for the regular-variation membership of $S_{n}$ and $M_{n}$. A common
shortcoming of formulas (\ref{Th1})--(\ref{Th4}) is the involved constants
which can not be accurately calculated in general. However, this is the
price we have to pay for highlighting the impact of the financial risk $%
Y_{i} $'s and weakening the moment conditions. Moreover, our explicit
expressions of $B_{n,i}$ and $D_{n,i}$ enable us to easily conduct numerical
estimates.

The following remarks and Corollary \ref{Cor1} contain some interesting
special cases of Theorem \ref{main2}, from which one can realize to some
extents the flexibility and generalization of our Theorem \ref{main2}.

\vskip0.2cm

\noindent \textbf{Remark 2.1. }If $\alpha =0$ then assertion (i) gives%
\begin{equation*}
\mathbb{P}\left( S_{n}>x\right) \sim \mathbb{P}\left( M_{n}>x\right) \sim 
\mathbb{P}\left( X_{n}\prod\limits_{j=1}^{n}Y_{j}>x\right) 
\end{equation*}%
and assertion (ii) reduces to%
\begin{equation*}
\mathbb{P}\left( S_{n}>x\right) \sim \mathbb{P}\left( M_{n}>x\right) \sim
\sum_{i=1}^{n}\mathbb{P}\left( X_{i}\prod\limits_{j=1}^{i}Y_{j}>x\right)
-\sum_{i=1}^{n-1}\mathbb{P}\left( \prod\limits_{j=1}^{i}Y_{j}>x\right) .
\end{equation*}

\vskip0.2cm

\noindent \textbf{Remark 2.2. }Clearly, if $\mathbb{E}\left\vert
X_{i}\right\vert ^{\beta }<\infty $ for every $i\geq 1$ and some $\beta
>\alpha $ then the two special conditions of assertion (i) hold in view of
Lemma \ref{Paul}(a) below. In this case, the last term of (\ref{Th1}) and (%
\ref{Th2}) can be expanded as follows by Breiman's lemma; see Breiman (1965),%
\begin{equation*}
\mathbb{P}\left( X_{n}\prod\limits_{j=1}^{n}Y_{j}>x\right) \sim \mathbb{E}%
\left( X_{n}\right) _{+}^{\alpha }\cdot \mathbb{P}\left(
\prod\limits_{j=1}^{n}Y_{j}>x\right) .
\end{equation*}%
Plugging this relation into (\ref{Th1}) and (\ref{Th2}) and noting that $%
\mathbb{E}\left( X_{n}\right) _{+}^{\alpha }=B_{n,n}=D_{n,n}$ yield%
\begin{equation*}
\mathbb{P}\left( S_{n}>x\right) \sim \sum_{i=1}^{n}B_{n,i}\mathbb{P}\left(
\prod\limits_{j=1}^{i}Y_{j}>x\right) \text{ and }\mathbb{P}\left(
M_{n}>x\right) \sim \sum_{i=1}^{n}D_{n,i}\mathbb{P}\left(
\prod\limits_{j=1}^{i}Y_{j}>x\right) .
\end{equation*}

\vskip0.2cm

\noindent \textbf{Remark 2.3. }By the proofs of Theorem \ref{main2}(i) and
Lemma \ref{convolution} below, if $X_{i}$'s are independent then (\ref%
{assume}) in assertion (i) can be weakened to $\overline{F}_{i}(x)=o\left( 
\overline{G}_{i+1}(x)\right) $.

\vskip0.2cm

In what follows, for a sequence $\{Z_{i};i\geq 1\}$ of iid rv's, we always
denote by $Z$ its generic rv.

\noindent \textbf{Remark 2.4. }By Lemma \ref{Paul}(a), if both $%
\{X_{i};i\geq 1\}$ and $\{Y_{i};i\geq 1\}$ are sequences of iid rv's then
only $\overline{F}(x)=o\left( \overline{G}(x)\right) $ suffices for
assertion (i). Moreover, we have%
\begin{equation*}
B_{n,i}=B_{n-i}=\mathbb{E}\left( X_{1}+S_{n-i}^{(2)}\right) _{+}^{\alpha }-%
\mathbb{E}\left( S_{n-i}^{(2)}\right) _{+}^{\alpha }=\mathbb{E}\left(
S_{n-i+1}\right) _{+}^{\alpha }\left( \mathbb{E}Y^{\alpha }\right) ^{-1}-%
\mathbb{E}\left( S_{n-i}\right) _{+}^{\alpha },
\end{equation*}%
and%
\begin{equation*}
D_{n,i}=D_{n-i}=\mathbb{E}\left( X_{1}+M_{n-i}^{(2)}\right) _{+}^{\alpha }-%
\mathbb{E}\left( M_{n-i}^{(2)}\right) ^{\alpha }=\mathbb{E}M_{n-i+1}^{\alpha
}\left( \mathbb{E}Y^{\alpha }\right) ^{-1}-\mathbb{E}M_{n-i}^{\alpha }.
\end{equation*}

\vskip0.2cm

\noindent \textbf{Remark 2.5. }The conditions of assertion (ii) do not
exclude the simultaneous occurrence of $\overline{F}_{i}(x)=o\left( 
\overline{G}_{i+1}(x)\right) $ for every $i\geq 1$. In such an
intersectional case, Lemma \ref{Paul}(b) and Remark 2.3 imply that assertion
(i) also holds and, hence, (\ref{Th3}) and (\ref{Th4}) should be equivalent
to (\ref{Th1}) and (\ref{Th2}), respectively. The latter fact can be easily
shown through Lemma \ref{small} below. Actually, for every $1\leq i\leq n-1$%
, by $\overline{F}_{i}(x)=o\left( \overline{G}_{i+1}(x)\right) $ and Lemma %
\ref{small}, we have%
\begin{equation*}
\mathbb{P}\left( X_{i}\prod\limits_{j=1}^{i}Y_{j}>x\right) -\mathbb{E}\left(
X_{i}\right) _{+}^{\alpha }\cdot \mathbb{P}\left(
\prod\limits_{j=1}^{i}Y_{j}>x\right) =o(1)\mathbb{P}\left(
\prod\limits_{j=1}^{i+1}Y_{j}>x\right) .
\end{equation*}%
On the other hand, it follows from Fatou's lemma that, for every $1\leq
i\leq n-1$,%
\begin{equation*}
\mathbb{P}\left( X_{n}\prod\limits_{j=1}^{n}Y_{j}>x\right) \gtrsim \mathbb{E}%
\left( X_{n}\prod\limits_{j=i+2}^{n}Y_{j}\right) _{+}^{\alpha }\cdot \mathbb{%
P}\left( \prod\limits_{j=1}^{i+1}Y_{j}>x\right) .
\end{equation*}%
Hence, 
\begin{equation*}
\sum_{i=1}^{n-1}\left( \mathbb{P}\left(
X_{i}\prod\limits_{j=1}^{i}Y_{j}>x\right) -\mathbb{E}\left( X_{i}\right)
_{+}^{\alpha }\cdot \mathbb{P}\left( \prod\limits_{j=1}^{i}Y_{j}>x\right)
\right) =o(1)\mathbb{P}\left( X_{n}\prod\limits_{j=1}^{n}Y_{j}>x\right) ,
\end{equation*}%
which implies that (\ref{Th3}) and (\ref{Th4}) are equivalent to (\ref{Th1})
and (\ref{Th2}), respectively.

\vskip0.2cm

The following corollary concerns another special case of Theorem \ref{main2}%
, in which the more explicit asymptotics can be derived. The assertion for $%
M_{n}$ was partially given by Theorem 6.1 of Tang and Tsitsiashvili (2003).
Recall that a real-valued rv $Z$ with survival function $\overline{H}$\ is
said to belong to the class $\mathcal{S}(\alpha )$ for some $\alpha \geq 0$
if%
\begin{equation}
\lim_{x\rightarrow \infty }\frac{\overline{H}(x-y)}{\overline{H}(x)}=\mathrm{%
e}^{\alpha y},\qquad y\in (-\infty ,\infty ),  \label{L}
\end{equation}%
and 
\begin{equation*}
\lim_{x\rightarrow \infty }\frac{\overline{H_{+}^{2\ast }}(x)}{\overline{H}%
(x)}=2\mathbb{E}\mathrm{e}^{\alpha Z}<\infty ,
\end{equation*}%
where $H_{+}(x)=H(x)\mathbf{1}_{\{x\geq 0\}}$ and $H_{+}^{2\ast }$ stands
for the $2$-fold convolution of $H_{+}$. In the literature, relation (\ref{L}%
) itself defines a larger class denoted by $\mathcal{L}(\alpha )$. See,
e.g., Cline (1987) and Pakes (2004, 2007) for more details on the classes $%
\mathcal{S}(\alpha )$ and $\mathcal{L}(\alpha )$. Note that, for a positive
rv $Z$, $\ln Z\in \mathcal{S}(\alpha )$ implies $Z\in \mathcal{R}_{-\alpha }$
and $\mathbb{E}Z^{\alpha }<\infty $.

\begin{corollary}
\label{Cor1}Under Assumption $\mathbf{A}$, let both $\{X_{i};i\geq 1\}$ and $%
\{Y_{i};i\geq 1\}$ be sequences of iid rv's. If $\ln Y\in \mathcal{S}(\alpha
)$ for some $\alpha \geq 0$ and $\lim_{x\rightarrow \infty }\overline{F}(x)/%
\overline{G}(x)=\theta \in \lbrack 0,\infty )$ then, for every $n\geq 1$,%
\begin{equation}
\mathbb{P}\left( S_{n}>x\right) \sim K_{n}\overline{G}(x)\text{ and }\mathbb{%
P}\left( M_{n}>x\right) \sim L_{n}\overline{G}(x),  \label{Co1}
\end{equation}%
where%
\begin{equation*}
K_{n}=\sum_{i=1}^{n}\left( \mathbb{E}\left( S_{n-i+1}\right) _{+}^{\alpha
}\left( \mathbb{E}Y^{\alpha }\right) ^{i-2}+\theta \left( \mathbb{E}%
Y^{\alpha }\right) ^{i}\right) \text{ and }L_{n}=\sum_{i=1}^{n}\left( 
\mathbb{E}M_{n-i+1}^{\alpha }\left( \mathbb{E}Y^{\alpha }\right)
^{i-2}+\theta \left( \mathbb{E}Y^{\alpha }\right) ^{i}\right) .
\end{equation*}%
Particularly, if $\alpha =0$ then, for every $n\geq 1$,%
\begin{equation*}
\mathbb{P}\left( S_{n}>x\right) \sim \mathbb{P}\left( M_{n}>x\right) \sim
\left( \theta +1\right) n\overline{G}(x).
\end{equation*}
\end{corollary}

\section{Lemmas and Proofs}

The following result is due to Corollary 2.1 of Hashorva and Li (2013),
which is motivated by Lemma 7.1 of Rootz\'{e}n (1986); see also Rootz\'{e}n
(1987). Note that for iid $Z_{i}$'s such that $\mathbb{P}(Z>x)\sim
cx^{-\alpha }$ the assertion was shown in Lemma 4.1(4) of Jessen and Mikosch
(2006).

\begin{lemma}
\label{product}Let $Z_{1},\ldots ,Z_{n}$ be $n$ positive and independent
rv's. If, for every $1\leq i\leq n$, $\mathbb{P}(Z_{i}>x)\sim \ell _{i}(\ln
x)(\ln x)^{\gamma _{i}-1}x^{-\alpha }$ for some positive constants $\alpha $%
, $\gamma _{i}$ and some slowly varying function $\ell _{i}(\cdot )$ then we
have 
\begin{equation*}
\mathbb{P}\left( \prod_{i=1}^{n}Z_{i}>x\right) \sim \frac{\alpha
^{n-1}\prod_{i=1}^{n}\Gamma (\gamma _{i})}{\Gamma \left(
\sum_{i=1}^{n}\gamma _{i}\right) }\left( \prod\limits_{i=1}^{n}\ell _{i}(\ln
x)\right) (\ln x)^{\sum_{i=1}^{n}\gamma _{i}-1}x^{-\alpha }.
\end{equation*}
\end{lemma}

\noindent \textbf{Proof of Theorem \ref{main3}}: The last relation in (\ref%
{Th5}) follows immediately from Lemma \ref{product}.\textit{\ }It remains to
verify that both the tails of $S_{n}$ and $M_{n}$ are asymptotically
equivalent to the right-hand side of (\ref{Th5}).\textit{\ }We only prove
the assertion for $S_{n}$, since the counterpart of $M_{n}$ can be obtained
similarly.

By Lemma \ref{product}, it is clear that the assertion holds for $%
S_{1}=X_{1}Y_{1}$. Now we assume by induction that the assertion holds for $%
n-1\geq 1$ and prove it for $n$. Recalling (\ref{S-M}), it holds that%
\begin{equation}
\mathbb{P}\left( S_{n}>x\right) =\mathbb{P}\left( Y_{1}\left(
X_{1}+S_{n-1}^{(2)}\right) >x\right) .  \label{m1}
\end{equation}%
From the induction assumption, we know that $S_{n-1}^{(2)}\in \mathcal{R}%
_{-\alpha }$ and $\overline{F}_{1}(x)=o(1)\mathbb{P}\left(
S_{n-1}^{(2)}>x\right) $. Noting also that $\overline{F}_{1}\in \mathcal{R}%
_{-\alpha }$ and $X_{1}$ is independent of $S_{n-1}^{(2)}$, we have (see,
e.g., Feller (1971), pp. 278) 
\begin{eqnarray*}
\mathbb{P}\left( X_{1}+S_{n-1}^{(2)}>x\right) &\sim &\mathbb{P}\left(
S_{n-1}^{(2)}>x\right) \\
&\sim &\frac{\alpha ^{n-1}\Gamma (\gamma ^{\ast })\prod_{i=2}^{n}\Gamma
(\gamma _{i})}{\Gamma \left( \gamma ^{\ast }+\sum_{i=2}^{n}\gamma
_{i}\right) }\ell _{n}^{\ast }(\ln x)\left( \prod\limits_{i=2}^{n}\ell
_{i}(\ln x)\right) \mathbb{\!}(\ln x)^{\gamma ^{\ast }+\sum_{i=2}^{n}\gamma
_{i}-1}x^{-\alpha }.
\end{eqnarray*}%
Then, applying Lemma \ref{product} to $Y_{1}$ and $X_{1}+S_{n-1}^{(2)}$ in (%
\ref{m1}) completes the proof. \hfill $\Box $

The next lemma is a restatement of the Corollary of Theorem 3 in Embrechts
and Goldie (1980).

\begin{lemma}
\label{Paul}Let $Y$ be a positive rv with survival function $\overline{G}\in 
\mathcal{R}_{-\alpha }$ for some $\alpha \geq 0$ and let $Z$ be a
real-valued rv with survival function $\overline{H}$. Assume that $Y$ and $Z$
are independent. Then $YZ\in \mathcal{R}_{-\alpha }$ if either (a) $%
\overline{H}(x)=o(\overline{G}(x))$ or (b) $\overline{H}\in \mathcal{R}%
_{-\alpha }$.
\end{lemma}

The first assertion of Lemma \ref{convolution} below is borrowed from Lemma
3.3 of Hao and Tang (2012); see also Lemma 4.4.2 of Samorodnitsky and Taqqu
(1994), and the second assertion is a special case of Proposition 2 of
Rogozin and Sgibnev (1999).

\begin{lemma}
\label{convolution}Let $Y$ and $Z$ be two real-valued rv's with survival
functions $\overline{G}$ and $\overline{H}$, respectively. If $\overline{G}%
\in \mathcal{R}_{-\alpha }$ for some $\alpha \geq 0$ and%
\begin{equation}
\mathbb{P}\left( \left\vert Z\right\vert >x\right) =o\left( \overline{G}%
(x)\right)  \label{assumption2}
\end{equation}%
then 
\begin{equation*}
\mathbb{P}\left( Y+Z>x\right) \sim \overline{G}(x).
\end{equation*}%
Particularly, if $Y$ and $Z$ are independent then (\ref{assumption2}) can be
weakened as $\overline{H}(x)=o\left( \overline{G}(x)\right) $.
\end{lemma}

Lemma \ref{key} below is crucial for the proof of our main theorem.

\begin{lemma}
\label{key}Let $Y$ be a positive rv with survival function $\overline{G}\in 
\mathcal{R}_{-\alpha }$ for some $\alpha \geq 0$ and let $Z_{1},\ldots
,Z_{n} $ be $n$\ real-valued rv's satisfying $\mathbb{E}\left( Z_{i}\right)
_{+}^{\alpha }<\infty $ for every $1\leq i\leq n$ and%
\begin{equation}
\mathbb{P}\left( \sum_{i=1}^{n}Z_{i}>x\right) \sim \sum_{i=1}^{n}c_{i}%
\mathbb{P}\left( Z_{i}>x\right)  \label{key1}
\end{equation}%
for $n$ nonnegative constants $c_{1},\ldots ,c_{n}$ such that $\max_{1\leq
i\leq n}c_{i}>0$. Assume further that $Y$ and $\{Z_{1},\ldots ,Z_{n}\}$ are
independent. Then%
\begin{equation}
\mathbb{P}\left( Y\sum_{i=1}^{n}Z_{i}>x\right) \sim \left( \mathbb{E}\left(
\sum_{i=1}^{n}Z_{i}\right) _{+}^{\alpha }-\sum_{i=1}^{n}c_{i}\mathbb{E}%
\left( Z_{i}\right) _{+}^{\alpha }\right) \mathbb{P}\left( Y>x\right)
+\sum_{i=1}^{n}c_{i}\mathbb{P}\left( YZ_{i}>x\right) .  \label{key2}
\end{equation}
\end{lemma}

One merit of Lemma \ref{key} is that we do not require $\mathbb{E}\left(
Z_{i}\right) _{+}^{\beta }<\infty $ for every $1\leq i\leq n$ and some $%
\beta >\alpha $. In return, the tails of products $\mathbb{P}\left(
YZ_{i}>x\right) $ for $1\leq i\leq n$ can not be expanded further.
Otherwise, relation (\ref{key2}) will reduce to Breiman's formula. If $Z_{i}$%
's are independent then relation (\ref{key1}) with $c_{1}=\cdots =c_{n}=1$
is usually called the max-sum equivalence property; see, e.g., Cai and Tang
(2004) for some heavy-tailed distribution classes satisfying such a
property. Moreover, even under some special dependence structures, including
the pairwise negative dependence and (quasi) asymptotic independence,
relation (\ref{key1}) still holds with $c_{1}=\cdots =c_{n}=1$ for $Z_{i}$'s
belonging to certain heavy-tailed distribution classes; see Chen and Yuen
(2009), Geluk and Tang (2009), and\ Tang (2008), among others.

\noindent \textbf{Proof of Lemma \ref{key}:} For every $0<\varepsilon <1$,
by relation (\ref{key1}), there is some $M>0$ such that the relations%
\begin{equation}
(1-\varepsilon )\sum_{i=1}^{n}c_{i}\mathbb{P}\left( Z_{i}>x\right) \leq 
\mathbb{P}\left( \sum_{i=1}^{n}Z_{i}>x\right) \leq (1+\varepsilon
)\sum_{i=1}^{n}c_{i}\mathbb{P}\left( Z_{i}>x\right)  \label{k1}
\end{equation}%
hold for all $x\geq M$. By this large $M$, we rewrite the left-hand side of (%
\ref{key2}) as%
\begin{eqnarray*}
\mathbb{P}\left( Y\sum_{i=1}^{n}Z_{i}>x\right) &=&\mathbb{P}\left(
Y\sum_{i=1}^{n}Z_{i}>x,Y>\frac{x}{M}\right) +\mathbb{P}\left(
Y\sum_{i=1}^{n}Z_{i}>x,Y\leq \frac{x}{M}\right) \\
&=&I_{1}(M,x)+I_{2}(M,x).
\end{eqnarray*}%
Applying Remark 4.1(a) below to $I_{1}(M,x)$, we have, for $M$ large enough,%
\begin{equation}
1-\varepsilon \leq \lim_{x\rightarrow \infty }\frac{I_{1}(M,x)}{\mathbb{E}%
\left( \sum_{i=1}^{n}Z_{i}\right) _{+}^{\alpha }\cdot \mathbb{P}\left(
Y>x\right) }\leq 1+\varepsilon .  \label{k2}
\end{equation}%
Consider $I_{2}(M,x)=\int_{0}^{x/M}\mathbb{P}\left(
\sum_{i=1}^{n}Z_{i}>x/y\right) G(\mathrm{d}y)$. It follows from (\ref{k1})
that%
\begin{equation}
(1-\varepsilon )J(M,x)\leq I_{2}(M,x)\leq (1+\varepsilon )J(M,x),  \label{k3}
\end{equation}%
where%
\begin{equation*}
J(M,x)=\sum_{i=1}^{n}c_{i}\mathbb{P}\left( YZ_{i}>x\right)
-\sum_{i=1}^{n}c_{i}\mathbb{P}\left( YZ_{i}>x,Y>\frac{x}{M}\right) .
\end{equation*}%
Using Remark 4.1(a) again to each summand of the second summation, we obtain
that,\ for $M$ large enough,%
\begin{equation}
1-\varepsilon \leq \lim_{x\rightarrow \infty }\frac{J(M,x)}{%
\sum_{i=1}^{n}c_{i}\mathbb{P}\left( YZ_{i}>x\right) -\sum_{i=1}^{n}c_{i}%
\mathbb{E}\left( Z_{i}\right) _{+}^{\alpha }\cdot \mathbb{P}\left(
Y>x\right) }\leq 1+\varepsilon .  \label{k4}
\end{equation}%
Combining (\ref{k2})--(\ref{k4}) and noting the arbitrariness of $%
\varepsilon $ complete the proof.\hfill $\Box $

\begin{lemma}
\label{small}Let $Y$ be a positive rv with survival function $\overline{G}%
\in \mathcal{R}_{-\alpha }$ for some $\alpha \geq 0$ and let $Z_{1},Z_{2}$
be $2$ real-valued rv's with distribution functions $H_{1},H_{2}$ satisfying 
$\overline{H}_{1}(x)=o\left( \overline{H}_{2}(x)\right) $ and $\mathbb{E}%
\left( Z_{2}\right) _{+}^{\alpha }<\infty $. Assume that $Y$ and $%
\{Z_{1},Z_{2}\}$ are independent. Then%
\begin{equation*}
\mathbb{P}\left( YZ_{1}>x\right) -\mathbb{E}\left( Z_{1}\right) _{+}^{\alpha
}\cdot \overline{G}(x)=o(1)\mathbb{P}\left( YZ_{2}>x\right) .
\end{equation*}
\end{lemma}

\proof For every $0<\varepsilon <1$, since $\overline{H}_{1}(x)=o\left( 
\overline{H}_{2}(x)\right) $, there is some $M$ such that for all $x\geq M$
the relation $\overline{H}_{1}(x)\leq \varepsilon \overline{H}_{2}(x)$
holds. Write%
\begin{equation*}
\mathbb{P}\left( YZ_{1}>x\right) =\mathbb{P}\left( YZ_{1}>x,Y>\frac{x}{M}%
\right) +\mathbb{P}\left( YZ_{1}>x,Y\leq \frac{x}{M}\right)
=I_{1}(M,x)+I_{2}(M,x).
\end{equation*}%
By Remark 4.1(a), choosing $M$ large enough, it holds that%
\begin{equation}
\lim_{x\rightarrow \infty }\frac{I_{1}(M,x)-\mathbb{E}\left( Z_{1}\right)
_{+}^{\alpha }\cdot \overline{G}(x)}{\mathbb{E}\left( Z_{1}\right)
_{+}^{\alpha }\cdot \overline{G}(x)}\leq \varepsilon .  \label{sm1}
\end{equation}%
For $I_{2}(M,x)$, by conditioning on $Y$ and noting that $\overline{H}%
_{1}(x)\leq \varepsilon \overline{H}_{2}(x)$ for $x\geq M$, we have%
\begin{equation}
I_{2}(M,x)\leq \varepsilon \mathbb{P}\left( YZ_{2}>x,Y\leq \frac{x}{M}%
\right) \leq \varepsilon \mathbb{P}\left( YZ_{2}>x\right) .  \label{sm2}
\end{equation}%
Moreover, Fatou's lemma gives%
\begin{equation}
\mathbb{P}\left( YZ_{2}>x\right) \gtrsim \mathbb{E}\left( Z_{2}\right)
_{+}^{\alpha }\cdot \overline{G}(x).  \label{sm3}
\end{equation}%
Therefore,%
\begin{eqnarray*}
&&\limsup_{x\rightarrow \infty }\frac{\mathbb{P}\left( YZ_{1}>x\right) -%
\mathbb{E}\left( Z_{1}\right) _{+}^{\alpha }\cdot \overline{G}(x)}{\mathbb{P}%
\left( YZ_{2}>x\right) } \\
&=&\limsup_{x\rightarrow \infty }\left( \frac{I_{1}(M,x)-\mathbb{E}\left(
Z_{1}\right) _{+}^{\alpha }\cdot \overline{G}(x)}{\mathbb{E}\left(
Z_{1}\right) _{+}^{\alpha }\cdot \overline{G}(x)}\cdot \frac{\mathbb{E}%
\left( Z_{1}\right) _{+}^{\alpha }\cdot \overline{G}(x)}{\mathbb{P}\left(
YZ_{2}>x\right) }+\frac{I_{2}(M,x)}{\mathbb{P}\left( YZ_{2}>x\right) }\right)
\\
&\leq &\varepsilon \left( \frac{\mathbb{E}\left( Z_{1}\right) _{+}^{\alpha }%
}{\mathbb{E}\left( Z_{2}\right) _{+}^{\alpha }}+1\right) ,
\end{eqnarray*}%
where in the last step we used (\ref{sm1}), (\ref{sm3}), and (\ref{sm2}) in
turn. Noting the arbitrariness of $\varepsilon $ completes the proof.%
\endproof

\noindent \textbf{Proof of Theorem \ref{main2}(i):} We only derive relation (%
\ref{Th1}) which implies $S_{n}\in \mathcal{R}_{-\alpha }$ by Lemma \ref%
{Paul}(b), then the assertions regarding $M_{n}$ follow from the similar
procedures with obvious modifications.

We proceed by the mathematical induction. Trivially, relation (\ref{Th1})
holds for $n=1$ with a by-product%
\begin{equation*}
\mathbb{P}\left( S_{1}>x\right) \gtrsim \mathbb{E}\left( X_{1}\right)
_{+}^{\alpha }\cdot \mathbb{P}\left( Y_{1}>x\right) .
\end{equation*}%
Assume by induction that relation (\ref{Th1}) holds for $n-1\geq 1$ with%
\begin{equation*}
\mathbb{P}\left( S_{n-1}>x\right) \gtrsim \mathbb{E}\left(
X_{1}+S_{n-2}^{(2)}\right) _{+}^{\alpha }\cdot \mathbb{P}\left(
Y_{1}>x\right) .
\end{equation*}%
Now we consider $S_{n}$ and recall that relation (\ref{m1}) holds. Applying
the induction assumption to $\{Y_{2},\ldots ,Y_{n}\}$ and $\{X_{2},\ldots
,X_{n}\}$ leads to%
\begin{equation}
\mathbb{P}\left( S_{n-1}^{(2)}>x\right) \gtrsim \mathbb{E}\left(
X_{2}+S_{n-2}^{(3)}\right) _{+}^{\alpha }\cdot \mathbb{P}\left(
Y_{2}>x\right) .  \label{m2}
\end{equation}%
Combining (\ref{m2}) with (\ref{assume}) gives%
\begin{equation*}
\mathbb{P}\left( \left\vert X_{1}\right\vert >x\right) =o\left( 1\right) 
\mathbb{P}\left( S_{n-1}^{(2)}>x\right) ,
\end{equation*}%
which together with Lemma \ref{convolution} implies%
\begin{equation*}
\mathbb{P}\left( X_{1}+S_{n-1}^{(2)}>x\right) \sim \mathbb{P}\left(
S_{n-1}^{(2)}>x\right) .
\end{equation*}%
Applying Lemma \ref{key} to (\ref{m1}) with $Y$, $Z_{1}$, $Z_{2}$ replaced
by $Y_{1}$, $X_{1}$, $S_{n-1}^{(2)}$, respectively, and $c_{1}=0$, $c_{2}=1$%
, we have%
\begin{eqnarray}
\mathbb{P}\left( S_{n}>x\right) &\sim &\left( \mathbb{E}\left(
X_{1}+S_{n-1}^{(2)}\right) _{+}^{\alpha }-\mathbb{E}\left(
S_{n-1}^{(2)}\right) _{+}^{\alpha }\right) \mathbb{P}\left( Y_{1}>x\right) +%
\mathbb{P}\left( Y_{1}S_{n-1}^{(2)}>x\right)  \notag \\
&=&B_{n,1}\mathbb{P}\left( Y_{1}>x\right) +\mathbb{P}\left( \widehat{S}%
_{n-1}^{(2)}>x\right) ,  \label{m3}
\end{eqnarray}%
where $\widehat{S}_{n-1}^{(2)}$ stands for $S_{n-1}^{(2)}$ with $Y_{2}$
replaced by $Y_{1}Y_{2}$. Clearly, $\{Y_{1}Y_{2},Y_{3},\ldots ,Y_{n}\}$ and $%
\{X_{2},\ldots ,X_{n}\}$ also satisfy all the conditions of assertion (i).
Thus, using the induction assumption to $\widehat{S}_{n-1}^{(2)}$ yields%
\begin{equation}
\mathbb{P}\left( \widehat{S}_{n-1}^{(2)}>x\right) \sim
\sum_{i=2}^{n-1}B_{n,i}\mathbb{P}\left( \prod\limits_{j=1}^{i}Y_{j}>x\right)
+\mathbb{P}\left( X_{n}\prod\limits_{j=1}^{n}Y_{j}>x\right) .  \label{m4}
\end{equation}%
A combination of (\ref{m3}) and (\ref{m4}) gives relation (\ref{Th1}).\hfill 
$\Box $

\noindent \textbf{Proof of Theorem \ref{main2}(ii):} Similarly as before, we
only derive relation (\ref{Th3}) by the mathematical induction. Trivially,
relation (\ref{Th3}) holds for $n=1$. Assume by induction that relation (\ref%
{Th3}) holds for $n-1\geq 1$, which implies $S_{n-1}^{(2)}\in \mathcal{R}%
_{-\alpha }$. Since $F_{1}\in \mathcal{R}_{-\alpha }$ and $X_{1}$ is
independent of $S_{n-1}^{(2)}$, it holds that%
\begin{equation*}
\mathbb{P}\left( X_{1}+S_{n-1}^{(2)}>x\right) \sim \mathbb{P}\left(
X_{1}>x\right) +\mathbb{P}\left( S_{n-1}^{(2)}>x\right) .
\end{equation*}%
Now, applying Lemma \ref{key} to (\ref{m1}) with $Y$, $Z_{1}$, $Z_{2}$
replaced by $Y_{1}$, $X_{1}$, $S_{n-1}^{(2)}$, respectively, and $%
c_{1}=c_{2}=1$, we have%
\begin{eqnarray}
&&\mathbb{P}\left( S_{n}>x\right)   \notag \\
&\sim &\!\!\left( \!\mathbb{E}\!\left( X_{1}+S_{n-1}^{(2)}\right)
_{+}^{\alpha }\!-\mathbb{E}\!\left( X_{1}\right) _{+}^{\alpha }\!-\mathbb{E}%
\!\left( S_{n-1}^{(2)}\right) _{+}^{\alpha }\right) \!\mathbb{P}\!\left(
Y_{1}>x\right) \!+\mathbb{P}\!\left( X_{1}Y_{1}>x\right) \!+\mathbb{P}%
\!\left( Y_{1}S_{n-1}^{(2)}>x\right)   \notag \\
&=&\!\!\left( B_{n,1}-\mathbb{E}\left( X_{1}\right) _{+}^{\alpha }\right) 
\mathbb{P}\left( Y_{1}>x\right) +\mathbb{P}\left( X_{1}Y_{1}>x\right) +%
\mathbb{P}\left( \widehat{S}_{n-1}^{(2)}>x\right) .  \label{m6}
\end{eqnarray}%
Since $\{Y_{1}Y_{2},Y_{3},\ldots ,Y_{n}\}$ and $\{X_{2},\ldots ,X_{n}\}$
also satisfy all the conditions of assertion (ii), using the induction
assumption on $\widehat{S}_{n-1}^{(2)}$ yields%
\begin{equation}
\mathbb{P}\left( \widehat{S}_{n-1}^{(2)}>x\right) \sim
\sum_{i=2}^{n-1}\left( B_{n,i}-\mathbb{E}\left( X_{i}\right) _{+}^{\alpha
}\right) \mathbb{P}\left( \prod\limits_{j=1}^{i}Y_{j}>x\right)
+\sum_{i=2}^{n}\mathbb{P}\left( X_{i}\prod\limits_{j=1}^{i}Y_{j}>x\right) .
\label{m7}
\end{equation}%
A combination of (\ref{m6}) and (\ref{m7}) gives relation (\ref{Th3}).\hfill 
$\Box $

\noindent \textbf{Proof of Corollary \ref{Cor1}}: Since $\ln Y\in \mathcal{S}%
(\alpha )$ and $\lim_{x\rightarrow \infty }\overline{F}(x)/\overline{G}%
(x)=\theta $, we can derive by Proposition 2 of Rogozin and Sgibnev (1999)
that, for every $i\geq 1$,%
\begin{equation}
\mathbb{P}\left( X_{i}\prod\limits_{j=1}^{i}Y_{j}>x\right) \sim \left( i%
\mathbb{E}X_{+}^{\alpha }+\theta \mathbb{E}Y^{\alpha }\right) \left( \mathbb{%
E}Y^{\alpha }\right) ^{i-1}\overline{G}(x),  \label{c1}
\end{equation}%
and, particularly,%
\begin{equation}
\mathbb{P}\left( \prod\limits_{j=1}^{i}Y_{j}>x\right) \sim i\left( \mathbb{E}%
Y^{\alpha }\right) ^{i-1}\overline{G}(x).  \label{c2}
\end{equation}%
If $\theta =0$, i.e., $\overline{F}(x)=o\left( \overline{G}(x)\right) $,
then Remark 2.4 indicates that Theorem \ref{main2}(i) holds. Plugging (\ref%
{c1}) and (\ref{c2}) into (\ref{Th1}) and (\ref{Th2}), and then rearranging
the constants with keeping in mind the two relations specified in Remark
2.4, we obtain the relations in (\ref{Co1}) with $\theta =0$. On the other
hand, if $\theta >0$ then Theorem \ref{main2}(ii) is valid. Plugging (\ref%
{c1}) and (\ref{c2}) into (\ref{Th3}) and (\ref{Th4}), and then rearranging
the constants, we complete the proof.\hfill $\Box $

\section{Appendix}

In this section, we derive some asymptotic results for the constant weighted
sums of partial products of $Y_{i}$'s with the uniformity of the constant
weights; see Theorem \ref{main1} below. We first prepare two important
lemmas.

\begin{lemma}
\label{uniform0}Let $Y$ be a positive rv with survival function $\overline{G}%
\in \mathcal{R}_{-\alpha }$ for some $\alpha \geq 0$ and let $\mathcal{Z}%
=\{Z\}$ be a set of positive rv's satisfying $\inf \mathcal{Z}>0$ and $%
\mathbb{E}\left( \sup \mathcal{Z}\right) ^{\alpha }<\infty $, where $\inf
/\sup \mathcal{Z}=\inf /\sup_{Z\in \mathcal{Z}}Z$. Assume that $Y$ and $%
\mathcal{Z}$ are independent. Then it holds uniformly for $Z\in \mathcal{Z}$
that%
\begin{equation}
\lim_{M\rightarrow \infty }\lim_{x\rightarrow \infty }\frac{\mathbb{P}\left(
YZ>x,Y>x/M\right) }{\mathbb{E}Z^{\alpha }\cdot \overline{G}(x)}=1.
\label{unif0}
\end{equation}
\end{lemma}

\proof For every $M>1>\delta >0$ and $x>0$, we have%
\begin{eqnarray*}
\mathbb{P}\left( YZ>x,Y>\frac{x}{M}\right) &=&\mathbb{P}\left( Y>\frac{x}{M}%
,Z>M\right) +\mathbb{P}\left( YZ>x,0<Z\leq \delta \right) +\mathbb{P}\left(
YZ>x,\delta <Z\leq M\right) \\
&=&I_{1}(M,x)+I_{2}(M,x)+I_{3}(M,x).
\end{eqnarray*}%
Since $Y$ and $\mathcal{Z}$ are independent, it holds that%
\begin{eqnarray}
\lim_{M\rightarrow \infty }\lim_{x\rightarrow \infty }\sup_{Z\in \mathcal{Z}}%
\frac{I_{1}(M,x)+I_{2}(M,x)}{\mathbb{E}Z^{\alpha }\cdot \overline{G}(x)}
&\leq &\lim_{M\rightarrow \infty }\lim_{x\rightarrow \infty }\sup_{Z\in 
\mathcal{Z}}\frac{\mathbb{P}\left( Z>M\right) \overline{G}\left( x/M\right) +%
\mathbb{P}\left( Z\leq \delta \right) \overline{G}(x/\delta )}{\mathbb{E}%
Z^{\alpha }\cdot \overline{G}(x)}  \notag \\
&\leq &\lim_{M\rightarrow \infty }\lim_{x\rightarrow \infty }\frac{\mathbb{P}%
\left( \sup \mathcal{Z}>M\right) \overline{G}\left( x/M\right) +\mathbb{P}%
\left( \inf \mathcal{Z}\leq \delta \right) \overline{G}(x/\delta )}{\mathbb{E%
}\left( \inf \mathcal{Z}\right) ^{\alpha }\cdot \overline{G}(x)}  \notag \\
&=&\lim_{M\rightarrow \infty }\frac{\mathbb{P}\left( \sup \mathcal{Z}%
>M\right) M^{\alpha }+\mathbb{P}\left( \inf \mathcal{Z}\leq \delta \right)
\delta ^{\alpha }}{\mathbb{E}\left( \inf \mathcal{Z}\right) ^{\alpha }} 
\notag \\
&\leq &\frac{\mathbb{P}\left( \inf \mathcal{Z}\leq \delta \right) }{\mathbb{E%
}\left( \inf \mathcal{Z}\right) ^{\alpha }},  \label{U01}
\end{eqnarray}%
where in the third and the fourth steps we used $G\in \mathcal{R}_{-\alpha }$
and $\mathbb{E}\left( \sup \mathcal{Z}\right) ^{\alpha }<\infty $,
respectively. For $I_{3}(M,x)$, we have%
\begin{eqnarray}
&&\lim_{M\rightarrow \infty }\lim_{x\rightarrow \infty }\sup_{Z\in \mathcal{Z%
}}\left\vert \frac{I_{3}(M,x)}{\mathbb{E}Z^{\alpha }\cdot \overline{G}(x)}%
-1\right\vert  \notag \\
&\leq &\lim_{M\rightarrow \infty }\lim_{x\rightarrow \infty }\sup_{Z\in 
\mathcal{Z}}\frac{\left\vert \int_{\delta }^{M}\left( \overline{G}(x/y)/%
\overline{G}(x)-y^{\alpha }\right) \mathbb{P}\left( Z\in dy\right)
\right\vert +\mathbb{E}Z^{\alpha }\mathbf{1}_{\{Z>M\}\cup \{Z\leq \delta \}}%
}{\mathbb{E}Z^{\alpha }}  \notag \\
&\leq &\lim_{M\rightarrow \infty }\lim_{x\rightarrow \infty }\frac{%
\sup_{\delta <y\leq M}\left\vert \overline{G}(x/y)/\overline{G}(x)-y^{\alpha
}\right\vert +\mathbb{E}\left( \sup \mathcal{Z}\right) ^{\alpha }\mathbf{1}%
_{\{\sup \mathcal{Z}>M\}}+\mathbb{P}\left( \inf \mathcal{Z}\leq \delta
\right) \delta ^{\alpha }}{\mathbb{E}\left( \inf \mathcal{Z}\right) ^{\alpha
}}  \notag \\
&\leq &\frac{\mathbb{P}\left( \inf \mathcal{Z}\leq \delta \right) }{\mathbb{E%
}\left( \inf \mathcal{Z}\right) ^{\alpha }},  \label{U02}
\end{eqnarray}%
where in the last step we used Theorem 1.5.2 of Bingham et al. (1987) to
neglect the first term of the numerator as $x\rightarrow \infty $. Combining
(\ref{U01}) with (\ref{U02}) and noting the arbitrariness of $\delta $
complete the proof.\endproof

\noindent \textbf{Remark 4.1. }Going along the same lines of the above proof
with corresponding modifications, we can obtain two variants of Lemma \ref%
{uniform0}: Let $Y$ be that in Lemma \ref{uniform0} and let $\mathcal{Z}$ be
a set of real-valued rv's independent of $Y$, then (a) relation (\ref{unif0}%
) with $\mathbb{E}Z^{\alpha }$ replaced by $\mathbb{E}Z_{+}^{\alpha }$,
denoted by (4.1$^{\prime }$), holds for every fixed $Z$ with $\mathbb{E}%
Z_{+}^{\alpha }<\infty $; (b) relation (4.1$^{\prime }$) holds uniformly for 
$Z\in \mathcal{Z}$ if $\alpha >0$ and $0<\mathbb{E}\left( \inf \mathcal{Z}%
\right) _{+}^{\alpha }\leq \mathbb{E}\left( \sup \mathcal{Z}\right)
_{+}^{\alpha }<\infty $.

\vskip0.2cm

Using Lemma \ref{uniform0} and the same idea as in the proof of Lemma \ref%
{key}, we have the following:

\begin{lemma}
\label{uniform}In addition to the other conditions of Lemma \ref{uniform0},
if $\mathbb{P}\left( Z>x-1\right) \sim \mathbb{P}(Z>x)$ holds uniformly for $%
Z\in \mathcal{Z}$ then it holds uniformly for $Z\in \mathcal{Z}$ that%
\begin{equation*}
\mathbb{P}\left( Y\left( 1+Z\right) >x\right) \sim \left[ \mathbb{E}\left(
1+Z\right) ^{\alpha }-\mathbb{E}Z^{\alpha }\right] \mathbb{P}\left(
Y>x\right) +\mathbb{P}\left( YZ>x\right) .
\end{equation*}
\end{lemma}

\begin{theorem}
\label{main1}Let $\{Y_{i};i\geq 1\}$ be a sequence of positive and
independent rv's with survival functions $\overline{G}_{i}\in \mathcal{R}%
_{-\alpha }$ for every $i\geq 1$ and some $\alpha \geq 0$. Assume that $%
\mathbb{E}Y_{i}^{\alpha }<\infty $ for every $i\geq 2$. Then, for every $%
n\geq 1$ and $0<a\leq b<\infty $, it holds uniformly for $(c_{1},\ldots
,c_{n})\in \lbrack a,b]^{n}$ that%
\begin{equation}
\mathbb{P}\left( \sum_{i=1}^{n}c_{i}\prod\limits_{j=1}^{i}Y_{j}>x\right)
\sim \sum_{i=1}^{n}A_{n,i}\mathbb{P}\left(
\prod\limits_{j=1}^{i}Y_{j}>x\right) ,  \label{theorem1}
\end{equation}%
where%
\begin{equation*}
A_{n,i}=\mathbb{E}\left(
\sum_{k=i}^{n}c_{k}\prod\limits_{j=i+1}^{k}Y_{j}\right) ^{\alpha }-\mathbb{E}%
\left( \sum_{k=i+1}^{n}c_{k}\prod\limits_{j=i+1}^{k}Y_{j}\right) ^{\alpha }.
\end{equation*}%
Particularly, if $\alpha =1$ then it holds uniformly for $(c_{1},\ldots
,c_{n})\in \lbrack a,b]^{n}$ that%
\begin{equation*}
\mathbb{P}\left( \sum_{i=1}^{n}c_{i}\prod\limits_{j=1}^{i}Y_{j}>x\right)
\sim \sum_{i=1}^{n}c_{i}\mathbb{P}\left(
\prod\limits_{j=1}^{i}Y_{j}>x\right) ,
\end{equation*}%
and if $\alpha =0$ then it holds uniformly for $(c_{1},\ldots ,c_{n})\in
\lbrack a,b]^{n}$ that%
\begin{equation*}
\mathbb{P}\left( \sum_{i=1}^{n}c_{i}\prod\limits_{j=1}^{i}Y_{j}>x\right)
\sim \mathbb{P}\left( \prod\limits_{j=1}^{n}Y_{j}>x\right) .
\end{equation*}
\end{theorem}

\proof We prove relation (\ref{theorem1}) by mathematical induction. For $%
n=1 $, by Theorem 1.5.2 of Bingham et al. (1987), it holds uniformly for $%
c_{1}\in \lbrack a,b]$ that%
\begin{equation*}
\mathbb{P}\left( c_{1}Y_{1}>x\right) \sim c_{1}^{\alpha }\mathbb{P}\left(
Y_{1}>x\right) =A_{1,1}\mathbb{P}\left( Y_{1}>x\right) .
\end{equation*}%
Hence, the assertion holds for $n=1$. Now we assume by induction that the
assertion holds for $n-1\geq 1$ and prove it for $n$. Define a set of
positive rv's as%
\begin{equation*}
\mathcal{Z}=\left\{ \sum_{i=2}^{n}\frac{c_{i}}{c_{1}}%
\prod_{j=2}^{i}Y_{j}:(c_{1},\ldots ,c_{n})\in \lbrack a,b]^{n}\right\} .
\end{equation*}%
It follows from Lemma \ref{Paul}(b) that $\prod_{j=2}^{i}Y_{j}\in \mathcal{R}%
_{-\alpha }\subset \mathcal{L}(0)$ for every $2\leq i\leq n$. Observing that 
$(c_{2}/c_{1},\ldots ,c_{n}/c_{1})\in \lbrack a/b,b/a]^{n-1}$, we obtain by
the induction assumption that, uniformly for $(c_{1},\ldots ,c_{n})\in
\lbrack a,b]^{n}$,%
\begin{eqnarray*}
\mathbb{P}\left( \sum_{i=2}^{n}\frac{c_{i}}{c_{1}}\prod_{j=2}^{i}Y_{j}>x-1%
\right) &\sim &\sum_{i=2}^{n}c_{1}^{-\alpha }A_{n,i}\mathbb{P}\left(
\prod\limits_{j=2}^{i}Y_{j}>x-1\right) \\
&\sim &\sum_{i=2}^{n}c_{1}^{-\alpha }A_{n,i}\mathbb{P}\left(
\prod\limits_{j=2}^{i}Y_{j}>x\right) \\
&\sim &\mathbb{P}\left( \sum_{i=2}^{n}\frac{c_{i}}{c_{1}}%
\prod_{j=2}^{i}Y_{j}>x\right) .
\end{eqnarray*}%
Moreover, it is obvious that%
\begin{equation*}
\inf \mathcal{Z}=\sum_{i=2}^{n}\frac{a}{b}\prod_{j=2}^{i}Y_{j}>0\text{ and }%
\mathbb{E}\left( \sup \mathcal{Z}\right) ^{\alpha }=\mathbb{E}\left(
\sum_{i=2}^{n}\frac{b}{a}\prod_{j=2}^{i}Y_{j}\right) ^{\alpha }<\infty .
\end{equation*}%
Hence, $\mathcal{Z}$ satisfies the conditions of Lemma \ref{uniform}, which
implies that, uniformly for $(c_{1},\ldots ,c_{n})\in \lbrack a,b]^{n}$,%
\begin{eqnarray}
\mathbb{P}\left( \sum_{i=1}^{n}c_{i}\prod\limits_{j=1}^{i}Y_{j}>x\right) &=&%
\mathbb{P}\left( Y_{1}\left( 1+\sum_{i=2}^{n}\frac{c_{i}}{c_{1}}%
\prod\limits_{j=2}^{i}Y_{j}\right) >\frac{x}{c_{1}}\right)  \notag \\
&\sim &c_{1}^{-\alpha }A_{n,1}\mathbb{P}\left( Y_{1}>\frac{x}{c_{1}}\right) +%
\mathbb{P}\left( Y_{1}\sum_{i=2}^{n}\frac{c_{i}}{c_{1}}\prod%
\limits_{j=2}^{i}Y_{j}>\frac{x}{c_{1}}\right)  \notag \\
&\sim &A_{n,1}\mathbb{P}\left( Y_{1}>x\right) +\mathbb{P}\left(
\sum_{i=2}^{n}c_{i}Y_{1}\prod\limits_{j=2}^{i}Y_{j}>x\right) .  \label{T1}
\end{eqnarray}%
For the second term of (\ref{T1}), regarding $Y_{1}Y_{2}$ as a whole and
using the induction assumption on $Y_{1}Y_{2},Y_{3},\ldots ,Y_{n}$, we have,
uniformly for $(c_{2},\ldots ,c_{n})\in \lbrack a,b]^{n-1}$,%
\begin{equation}
\mathbb{P}\left(
\sum_{i=2}^{n}c_{i}Y_{1}\prod\limits_{j=2}^{i}Y_{j}>x\right) \sim
\sum_{i=2}^{n}A_{n,i}\mathbb{P}\left( \prod\limits_{j=1}^{i}Y_{j}>x\right) .
\label{T2}
\end{equation}%
A combination of (\ref{T1}) and (\ref{T2}) completes the proof.\endproof

Similarly as in Corollary \ref{Cor1}, assuming further that $\{Y_{i};i\geq
1\}$ is a sequence of iid rv's and $\ln Y\in \mathcal{S}(\alpha )$ for some $%
\alpha \geq 0$ leads to a series of explicit results. We conclude them in
the following Corollary \ref{Cor2}.

\begin{corollary}
\label{Cor2}Let $\{Y_{i};i\geq 1\}$ be a sequence of positive and iid rv's
with common survival function $\overline{G}$. If $\ln Y\in \mathcal{S}%
(\alpha )$ for some $\alpha \geq 0$ then, for every $n\geq 1$ and $0<a\leq
b<\infty $, it holds uniformly for $(c_{1},\ldots ,c_{n})\in \lbrack
a,b]^{n} $ that%
\begin{equation*}
\mathbb{P}\left( \sum_{i=1}^{n}c_{i}\prod\limits_{j=1}^{i}Y_{j}>x\right)
\sim \sum_{i=1}^{n}\mathbb{E}\left(
\sum_{k=i}^{n}c_{k}\prod\limits_{j=1}^{k-i+1}Y_{j}\right) ^{\alpha }\left( 
\mathbb{E}Y^{\alpha }\right) ^{i-2}\cdot \overline{G}(x).
\end{equation*}%
Particularly, if $\alpha =1$ then it holds uniformly for $(c_{1},\ldots
,c_{n})\in \lbrack a,b]^{n}$ that%
\begin{equation*}
\mathbb{P}\left( \sum_{i=1}^{n}c_{i}\prod\limits_{j=1}^{i}Y_{j}>x\right)
\sim \sum_{i=1}^{n}ic_{i}\left( \mathbb{E}Y\right) ^{i-1}\cdot \overline{G}%
(x),
\end{equation*}%
and if $\alpha =0$ then it holds uniformly for $(c_{1},\ldots ,c_{n})\in
\lbrack a,b]^{n}$ that%
\begin{equation*}
\mathbb{P}\left( \sum_{i=1}^{n}c_{i}\prod\limits_{j=1}^{i}Y_{j}>x\right)
\sim n\overline{G}(x).
\end{equation*}
\end{corollary}

\vskip0.5cm

\noindent \textbf{Acknowledgments.} The authors are very grateful to an
anonymous reviewer for her/his very thorough reading of the paper and
valuable suggestions. E. Hashorva kindly acknowledges partial support from
the Swiss National Science Foundation Project 200021-140633/1 and from the
project RARE -318984 (a Marie Curie International Research Staff Exchange
Scheme Fellowship within the 7th European Community Framework Programme). J.
Li kindly acknowledges the financial support from the National Natural
Science Foundation of China (Grant No. 11201245) and partial support from
Swiss National Science Foundation Grant 200021-134785.

\end{document}